\documentclass[preprint,11pt]{elsarticle}
\usepackage{graphicx}
\usepackage[usenames,dvipsnames]{color}

\begin{document}
\newcommand{\up}{\vspace*{-0.05cm}}
\newcommand{\qqed}{\hfill$\rule{.05in}{.1in}$\vspace{.3cm}}
\newcommand{\pf}{\noindent{\bf Proof: }}
\newtheorem{thm}{Theorem}
\newtheorem{lem}{Lemma}
\newtheorem{prop}{Proposition}
\newtheorem{prob}{Problem}
\newtheorem{ex}{Example}
\newtheorem{cor}{Corollary}
\newtheorem{conj}{Conjecture}
\newtheorem{cl}{Claim}
\newtheorem{df}{Definition}
\newtheorem{rem}{Remark}
\newcommand{\beq}{\begin{equation}}
\newcommand{\eeq}{\end{equation}}
\newcommand{\<}[1]{\left\langle{#1}\right\rangle}
\newcommand{\be}{\beta}
\newcommand{\ee}{\end{enumerate}}
\newcommand{\Bul}{\mbox{$\bullet$ } }
\newcommand{\al}{\alpha}
\newcommand{\ep}{\epsilon}
\newcommand{\si}{\sigma}
\newcommand{\om}{\omega}
\newcommand{\la}{\lambda}
\newcommand{\La}{\Lambda}
\newcommand{\Ga}{\Gamma}
\newcommand{\ga}{\gamma}
\newcommand{\im}{\Rightarrow}
\newcommand{\2}{\vspace{.2cm}}
\newcommand{\es}{\emptyset}

\journal{arXiv.org}

\begin{frontmatter}

\title{Braess' Paradox in a Generalised Traffic Network}
\date{}
\author[label2]{Vadim Zverovich\corref{cor1}}
\ead{vadim.zverovich@uwe.ac.uk}
\cortext[cor1]{Corresponding author}
\author[label2a]{Erel Avineri}
\address[label2]{Department of Mathematics and Statistics, University of the West of England, Bristol, BS16 1QY, UK}
\address[label2a]{Centre for Transport and Society, University of the West of England, Bristol, BS16 1QY, UK}

\begin{abstract}
The classical network configuration introduced by Braess in 1968
is of fundamental significance because Valiant and Roughgarden
showed in 2006 that `the ``global" behaviour of an equilibrium
flow in a large random network is similar to that in Braess'
original four-node example'. In this paper, a natural
generalisation of Braess' network is introduced and conditions for
the occurrence of Braess' paradox are formulated for the
generalised network.

The Braess' paradox has been studied mainly in the context of the
classical problem introduced by Braess and his colleagues,
assuming a certain type of networks. Specifically, two pairs of
links in those networks are assumed to have the same volume-delay
functions. The occurrence of Braess' paradox for this specific
case of network symmetry was investigated by Pas and Principio in
1997.  Such a symmetry is not common in real-life networks because
the parameters of volume-delay functions are associated with roads
physical and functional characteristics, which typically differ
from one link to another (e.g. roads in networks are of different
length).  Our research provides an extension of previous studies
on Braess' paradox by considering arbitrary volume-delay
functions, i.e. symmetry properties are not assumed for any of the
network's links and the occurrence of Braess' paradox is studied
for a general configuration.

\end{abstract}
\end{frontmatter}

\section{Introduction}
The traffic network can be generally described as a game,
where a finite number of interdependent network users compete with
each other by making simultaneous decisions route choices. It is
commonly assumed that network users non-cooperatively interact
with each other in the traffic network in order to minimize their
travel costs. This problem can be modelled as an $N$-person
nonzero-sum game (see \cite{Byu}). Its solution assumes the
existence of equilibrium. The concept of equilibrium in the
context of transport systems had appeared in the 1950s
(\cite{Bec,War}) and is based on the general assumption that
network users are making adjustments to their travel choices until
a state of equilibrium is reached, i.e when no individual can make
a further improvement as a result of any individual choice. A
specific situation investigated in this work is the effect of
introducing a new link to a congested traffic network, and the
likelihood of this additional capacity to improve the system's
performance (measured by users' aggregated cost or travel time).
The addition of a link to an existing transport system may lead to
undesired situations.

The well-known Braess' paradox \cite{Bra} describes situations
when adding a new link to a transport network might not reduce
congestion in the network but instead increase it. This is due to
individual entities acting selfishly/separately when making their
travel plan choices and hence forcing the system as a whole not to
operate optimally. Deeper insight into this paradox from the
viewpoint of the structure and characteristics of networks may
help transport planners to avoid the occurrence of Braess-like
situations in real-life networks. The same paradox also applies to
situations when an existing link is removed from a network. The
famous example is the closing of the 42nd street in New York City
in 1990. Contrary to all expectations, the closing reduced the
amount of congestion in the area \cite{Kol}.

Pas and Principio \cite{Pas} investigated the existence of Braess'
paradox in his classical network configuration, where travel times
of links are specified in such a way that the resulting network is
symmetrical. Their result describes the situations when Braess'
paradox occurs in the above-mentioned network but is only limited
to the type of networks shown in Figure 4, when the travel times
of some links are symmetrical and there are further restrictions
on free flow travel times and delay parameters (see Section 5).
However, there is no symmetry in typical real-life networks.

Steinberg and Zangwill \cite{Ste} showed that, under reasonable
assumptions, Braess' paradox is about as likely to occur as not
occur in general transportation networks. Under other assumptions,
Valiant and Roughgarden \cite{Val} proved that Braess' paradox is
likely to occur in a natural random network model. More precisely,
for a given appropriate total flow, they showed that in almost all
networks there is a set of links whose removal improves the travel
time at equilibrium, i.e. Braess' paradox is widespread. Thus, the
fundamental significance of Steinberg-Zangwill's and
Valiant-Roughgarden's results is that the presence of Braess'
paradox is not rare and exceptional, but rather widespread.
Therefore, real-life networks must be analysed from the viewpoint
of this phenomenon before adding/constructing a new link/road.
Moreover, Valiant and Roughgarden \cite{Val} showed that `the
``global" behaviour of an equilibrium flow in a large random
network is similar to that in Braess' original four-node example'.
This important result means that the classical network
configuration introduced by Braess in 1968 is fundamental and it
should be analysed in depth.

In this paper, a natural generalisation of Braess' network is
introduced in Section 2. Such a generalised network can be reduced
to the classical network configuration $N/N^+$, where $N^+$ is the
network of Figure 1 and the network $N$ is $N^+$ with the link
$(b,c)$ removed. However, in contrast to the classical case with
symmetrical travel times of links, the network configuration
$N/N^+$ has arbitrary linear travel time functions. The travel
times at equilibria in the networks $N$ and $N^+$ are completely
described in Section 3, while Section 4 provides necessary and
sufficient conditions for the occurrence of Braess' paradox in
these networks.

In Section 5, a known result describing the existence of Braess'
paradox in the symmetrical network configuration is obtained and
generalised. The motivation for such an extension was given by Pas
and Principio in their conclusions in \cite{Pas}. Also, the
concept of a pseudo-paradox is introduced and the conditions for
its occurrence in $N/N^+$ are given. In particular, it is shown
that in an asymmetrical network configuration the pseudo-paradox
is happening for any total flow. The proof of all the theorems is
given in Section \ref{Proof}.

\section{A Generalisation of Braess' Network Configuration and Notation}

The classical network configuration introduced by Braess
\cite{Bra} consists of three paths:
$$
P_1=a-b-d, \;\; P_2=a-c-d, \;\; P_3=a-b-c-d.
$$
This network is denoted by $N^+$, and it has four nodes and five
links, where $a$ is the origin of all travel demand, and $d$ is
the destination of all demand (see Figure 1). The network $N$ is
$N^+$ with the link $(b,c)$ removed.

In 2006 Valiant and Roughgarden \cite{Val} proved an important and
interesting result that Braess' paradox is likely to happen in a
natural random network model. Also, they showed that `the
``global" behaviour of an equilibrium flow in a large random
network is similar to that in Braess' original four-node example'.
Thus, Braess' network configuration is of fundamental
significance.

\unitlength=1.00mm
\begin{picture}(150,40)
\put(50,20){\circle*{2}} \put(70,30){\circle*{2}}
\put(70,10){\circle*{2}} \put(90,20){\circle*{2}}

\put(45,19){\makebox{$a$}} \put(93,19){\makebox{$d$}}
\put(69,32){\makebox{$b$}}\put(69,6){\makebox{$c$}}

\qbezier (50,20)(70,30)(70,30) \qbezier
(50,20)(70,10)(70,10)\qbezier (70,30)(70,10)(70,10)\qbezier
(70,30)(90,20)(90,20)\qbezier (70,10)(90,20)(90,20)


\qbezier (68,28)(69,29.5)(69,29.5)
\qbezier(67.3,29.5)(69,29.5)(69,29.5)

\qbezier (67.3,10.5)(69,10.5)(69,10.5)
\qbezier(68,12)(69,10.5)(69,10.5)

\qbezier (88,18)(89,19.5)(89,19.5)
\qbezier(87.3,19.5)(89,19.5)(89,19.5)

\qbezier (87.3,20.5)(89,20.5)(89,20.5)
\qbezier(88,22)(89,20.5)(89,20.5)

\qbezier(70.8,12.5)(70,11)(70,11)
\qbezier(69.2,12.5)(70,11)(70,11)


\end{picture}

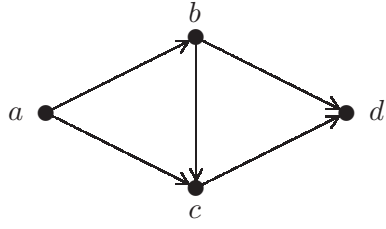
\vspace*{-.3cm} \centerline{{\bf Figure 1.} Braess' Network
Configuration.}

\bigskip

Let us consider a natural generalisation of Braess' network, where
every link of the network in Figure 1 is replaced by a path of
arbitrary length (i.e. a route with any number of links). Thus,
the generalised network `comprises' five paths of arbitrary
length: $(a,b)$-path, $(b,d)$-path, $(a,c)$-path, $(c,d)$-path and
$(b,c)$-path. This is illustrated by the first network in Figure
2, where the $(b,c)$-path is of length 3 and other four paths have
length 2. Every link $(i,j)$ in the resulting network has a linear
travel time function
$$
\al_{ij}+\be_{ij}f_{ij},
$$
where $\al_{ij}\ge 0$ is the free flow travel time for the link
$(i,j)$, $\be_{ij}>0$ is the delay parameter for $(i,j)$, and
$f_{ij}\ge 0$ is the flow on the link $(i,j)$. A fixed traffic
coming from outside the network is allowed. For example, in Figure
2 the dashed arrows represent a fixed traffic $\tilde{f}$ coming
to the network and then going outside.

The assumption of a linear relationship between traffic volume on
a link and the travel time on it (so-called `volume-delay
function') is common in the context of Braess' paradox. Although
there is evidence to support such a linear
approximation\cite{Wal}, many different types of volume-delay
functions, e.g. BPR functions \cite{Bur}, have been applied (see
\cite{Bran} for a review article). However, the investigation of
the non-linear case is not in the scope of this work.


\unitlength=1.00mm
\begin{picture}(200,40)

\put(10,20){\circle*{2}} \put(30,30){\circle*{2}}
\put(30,10){\circle*{2}} \put(50,20){\circle*{2}}

\put(20,25){\circle*{2}} \put(40,25){\circle*{2}}
\put(20,15){\circle*{2}} \put(40,15){\circle*{2}}

\put(30,17){\circle*{2}}\put(30,24){\circle*{2}}

\put(5,19){\makebox{$a$}} \put(53,19){\makebox{$d$}}
\put(29,32){\makebox{$b$}}\put(29,6){\makebox{$c$}}

\qbezier (10,20)(30,30)(30,30) \qbezier
(10,20)(30,10)(30,10)\qbezier (30,30)(30,10)(30,10)\qbezier
(30,30)(50,20)(50,20)\qbezier (30,10)(50,20)(50,20)

\qbezier (28,28)(29,29.5)(29,29.5)
\qbezier(27.3,29.5)(29,29.5)(29,29.5)

\qbezier (18,23)(19,24.5)(19,24.5)
\qbezier(17.3,24.5)(19,24.5)(19,24.5)

\qbezier (27.3,10.5)(29,10.5)(29,10.5)
\qbezier(28,12)(29,10.5)(29,10.5)

\qbezier (17.3,15.5)(19,15.5)(19,15.5)
\qbezier(18,17)(19,15.5)(19,15.5)

\qbezier (48,18)(49,19.5)(49,19.5)
\qbezier(47.3,19.5)(49,19.5)(49,19.5)

\qbezier (38,13)(39,14.5)(39,14.5)
\qbezier(37.3,14.5)(39,14.5)(39,14.5)

\qbezier (47.3,20.5)(49,20.5)(49,20.5)
\qbezier(48,22)(49,20.5)(49,20.5)

\qbezier (37.3,25.5)(39,25.5)(39,25.5)
\qbezier(38,27)(39,25.5)(39,25.5)

\qbezier(30.8,12.5)(30,11)(30,11)
\qbezier(29.2,12.5)(30,11)(30,11)

\qbezier(30.8,19.5)(30,18)(30,18)
\qbezier(29.2,19.5)(30,18)(30,18)

\qbezier(30.8,26.5)(30,25)(30,25)
\qbezier(29.2,26.5)(30,25)(30,25)


\qbezier(20,25)(20,28)(20,28)
\qbezier(20,29)(20,31)(20,31)\qbezier(20,32)(20,34)(20,34)

\qbezier(20.8,27.5)(20,26)(20,26)
\qbezier(19.2,27.5)(20,26)(20,26)

\put(17,29){\makebox{$\tilde{f}$}}

\qbezier(40,25)(40,28)(40,28)
\qbezier(40,29)(40,31)(40,31)\qbezier(40,32)(40,34)(40,34)

\qbezier(40.8,32.5)(40,34)(40,34)
\qbezier(39.2,32.5)(40,34)(40,34)

\put(41,29){\makebox{$\tilde{f}$}}

\put(68,18){\makebox{\huge $\Rightarrow$}}


\put(90,20){\circle*{2}} \put(110,30){\circle*{2}}
\put(110,10){\circle*{2}} \put(130,20){\circle*{2}}

\put(85,19){\makebox{$a$}} \put(133,19){\makebox{$d$}}
\put(109,32){\makebox{$b$}}\put(109,6){\makebox{$c$}}

\put(94,27){\makebox{\small $\al_1,\be_1$}}
\put(117,27){\makebox{\small $\al_2,\be_2$}}
\put(111,19){\makebox{\small $\al_3,\be_3$}}
\put(94,10.5){\makebox{\small $\al_4,\be_4$}}
\put(117,10.5){\makebox{\small $\al_5,\be_5$}}

\qbezier (90,20)(110,30)(110,30) \qbezier
(90,20)(110,10)(110,10)\qbezier (110,30)(110,10)(110,10)\qbezier
(110,30)(130,20)(130,20)\qbezier (110,10)(130,20)(130,20)

\qbezier (108,28)(109,29.5)(109,29.5)
\qbezier(107.3,29.5)(109,29.5)(109,29.5)

\qbezier (107.3,10.5)(109,10.5)(109,10.5)
\qbezier(108,12)(109,10.5)(109,10.5)

\qbezier (128,18)(129,19.5)(129,19.5)
\qbezier(127.3,19.5)(129,19.5)(129,19.5)

\qbezier (127.3,20.5)(129,20.5)(129,20.5)
\qbezier(128,22)(129,20.5)(129,20.5)

\qbezier(110.8,12.5)(110,11)(110,11)
\qbezier(109.2,12.5)(110,11)(110,11)

\end{picture}

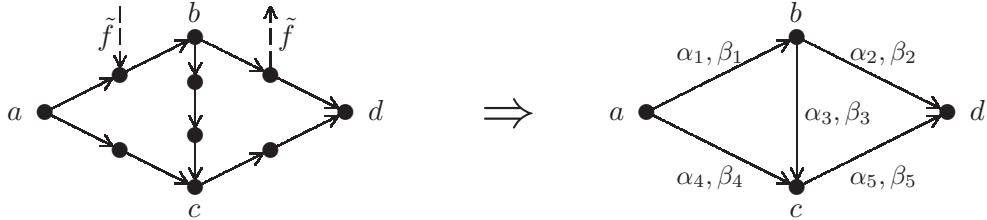
\vspace*{-.3cm} \centerline{{\bf Figure 2.} A generalised network
reduced to a four-node network $N^+$.}

\bigskip

Suppose we want to decide whether Braess' paradox occurs when
removing a link on the path going from $b$ to $c$. If a particular
link $(i,j)$ has a fixed  flow $\tilde{f}$ coming from outside the
network, i.e. not the internal flow $f$ going from $a$ to $d$
through this link, then its travel time function can be written as
follows:
$$
\al_{ij}+\be_{ij}f_{ij} = \al_{ij}+\be_{ij}(f+\tilde{f}) =
(\al_{ij}+\be_{ij}\tilde{f}) + \be_{ij}f = \tilde{\al}_{ij} +
\be_{ij}f.
$$
The updated function depends on the internal flow and it is
linear, because the external flow $\tilde{f}$ is fixed and hence
$\tilde{\al}_{ij}$ is a fixed number. Thus, the first step is to
update all travel time functions taking into accounts external
flows. Further, it is easy to see that the total travel time
functions for the above paths are linear functions, since the
travel time functions for links are linear and all the links
belonging to one of the paths share the same internal flow. For
instance, if such a path $P$ with an internal flow $f$ consists of
two links $(i,j)$ and $(j,k)$, then  the total travel time
function is as follows:
$$
\al_{ij}+\be_{ij}f + \al_{jk}+\be_{jk}f =
\al_{ij}+\al_{jk}+(\be_{ij}+\be_{jk})f = \al_{P}+\be_{P}f.
$$
Thus, if all the above paths are replaced by single links, then we
obtain Braess' network with arbitrary linear travel time functions
(see Figure 2):

$$\al_1+\be_1f_{ab}\;\;  \mbox{for link} \;(a,b);$$
$$\al_2+\be_2f_{bd}\;\;  \mbox{for link} \;(b,d);$$
$$\al_3+\be_3f_{bc}\;\;  \mbox{for link} \;(b,c);$$
$$\al_4+\be_4f_{ac}\;\;  \mbox{for link} \;(a,c);$$
$$\al_5+\be_5f_{cd}\;\;  \mbox{for link} \;(c,d).$$


Note that in Braess' original example \cite{Bra} and in many of
the studies that followed it (e.g. \cite{Byu}), the network is
symmetric, i.e the time functions for the links (a,b) and (c,d)
are same as well as the time functions for the links (b,d) and
(a,c), and the free flow travel times for the links $(a,b)$ and
$(c,d)$ are equal to zero. The occurrence of Braess' paradox in
this symmetrical network configuration was described by Pas and
Principio \cite{Pas} (see Corollary \ref{PP}).

In this paper, we consider a more general situation with arbitrary
linear time functions. The existence of Braess' paradox in such a
network can be decided by using Theorems \ref{T1}--\ref{T4}, where
the network $N$ is $N^+$ with the link $(b,c)$ removed.

Let $Q>0$ denote the total flow in $N/N^+$, i.e
$Q=f_{ab}+f_{ac}=f_{bd}+f_{cd}$. Note that $f_{ij}$ and $Q$ are
not necessarily integer numbers. Let us denote
$$\al_{ij}=\al_i+\al_j,$$ e.g. $\al_{12}$ means $\al_1+\al_2$, and
$\be_{ij}$ is defined similarly. Also,
$$
\al=\al_{45}-\al_{12}, \qquad \bar\al=\al_{4}-\al_{13}, \qquad
\hat\al=\al_{2}-\al_{35}, \qquad
$$
and
$$
\be = \be_{1245} = \be_1+\be_2+\be_4+\be_5.
$$
The following equality will be used throughout the paper:
$$
\al = \bar\al - \hat\al.
$$

Further, we introduce the Braess numbers ${\cal B}_i$ for
$i=1,2,3,4$:
$$
{\cal B}_1 = \be_1\be_5-\be_2\be_4, \quad {\cal B}_2 =
\be_{135}\be-\be_{12}\be_{45}, \quad
{\cal B}_3 = \be^2_{45}\be_{134}-\be^2_4\be \quad \mbox{and} \quad
{\cal B}_4 = \be^2_{12}\be_{235}-\be^2_2\be.
$$

Also, two parameters $\mu_1$ and $\mu_2$ are defined as follows:
$$
\mu_1 = {\hat{\al}\be_{14} - \al\be_{3} \over
\be_3\be_{45}+\be_5\be_{14}}, \qquad \mu_2 = {\bar{\al}\be_{25} +
\al\be_{3} \over \be_1\be_{25}+\be_3\be_{12}}.
$$

\section{Equilibria in $N$ and $N^+$}

It is well known that a user equilibrium always exists, and in a
network without capacities, it is essentially unique (e.g. see
\cite{Sch}).

\2

A network with one source and one sink is said to be at
equilibrium if
\begin{itemize}
    \item [(a)] The travel time on paths with non-vanishing flow is the same,
and it is denoted by $T_{Eq}$, and
    \item [(b)] The travel time on paths with no flow is at least $T_{Eq}$.
\end{itemize}

The equilibrium described above is associated with aggregated
strategic behaviour of all road users, described as an $N$-person
Wardrop/Nash equilibrium. In equilibrium, no user can decrease
his/her route travel time by unilaterally switching routes \cite
{War}. In other words, if a network is not at equilibrium, then
some users of the network (e.g. drivers) can switch their routes
in order to improve their travel time. However, if a driver
decides to switch to a better route, then the travel time for this
route increases, and, after a certain period of time, it will
become impossible to improve drivers' travel times by switching
the routes. Thus, the equilibrium describes `stable state'
behaviour in a network, and no driver has any incentive to switch
routes at equilibrium because it will not improve their current
travel times.

The following lemma describes the equilibrium in the network $N$,
which is $N^+$ with the link $(b,c)$ removed. Note that in Lemma
\ref{L1} the case (a) corresponds to the situation when the path
$P_1$ has a vanishing flow and $P_2$ has a non-vanishing flow in
$N$. In case (b) the path $P_1$ has a non-vanishing flow and $P_2$
has a vanishing flow, and in case (c) no path has a vanishing
flow. Also, the cases (a) and (b) in this lemma are mutually
exclusive because one of the numbers $-\al/\be_{45}$ and
$\al/\be_{12}$  is negative, or they both are equal to zero.

\begin{lem} \label{L1}
In the network $N$, the travel time at equilibrium is as follows:

\begin{itemize}

    \item [(a)] $T_{Eq} = \al_{45}+Q\be_{45}\;\;$ if $\;\; 0<Q\le
    {-\al/\be_{45}}$;

    \item [(b)] $T_{Eq} = \al_{12}+Q\be_{12}\;\;$ if $\;\; 0<Q\le { \al/\be_{12}}$;

    \item [(c)] $T_{Eq} = \al_{12}+(\al+Q\be_{45})\be_{12}/\be\;\;$ if
$\;\; Q>\max\{{\al/\be_{12}}; {-\al/\be_{45}}\}$.

\end{itemize}
\end{lem}

\begin{pf}
Let us denote $f_{ab}$ by $h$. Then $f_{bd}=h$ and
$f_{ac}=f_{cd}=Q-h$. We have
$$
T_1= \al_{12}+h\be_{12} \quad \mbox{and} \quad T_2=
\al_{45}+(Q-h)\be_{45},
$$
where $T_i$ is the travel time on the path $P_i$.

\2\noindent {\bf\em Case (a):} Suppose that $P_1$ has a vanishing
flow and $P_2$ has a non-vanishing flow. Then $Q>h=0$ and
$$
T_1= \al_{12} \quad \mbox{and} \quad T_2= \al_{45}+Q\be_{45}.
$$
At equilibrium, $T_1\ge T_2$, i.e. $\al_{12} \ge
\al_{45}+Q\be_{45}$ or $Q\le -\al/\be_{45}$. The travel time at
equilibrium is
$$
T_{Eq} = T_2 = \al_{45}+Q\be_{45}.
$$

\2 \noindent {\bf\em Case (b):} Assume that $P_1$ has a
non-vanishing flow and $P_2$ has a vanishing flow. We have $Q=h>0$
and
$$
T_1= \al_{12}+Q\be_{12} \quad \mbox{and} \quad T_2= \al_{45}.
$$
At equilibrium, $T_1\le T_2$, i.e. $\al_{12}+Q\be_{12} \le
\al_{45}$ or $Q\le \al/\be_{12}$. The travel time at equilibrium
is as follows:
$$
T_{Eq} = T_1 = \al_{12}+Q\be_{12}.
$$

\2\noindent {\bf\em Case (c):} Suppose that no path has a
vanishing flow. We have $Q>h>0$. At equilibrium, $T_1 = T_2$, i.e.
$$
\al_{12}+h\be_{12} = \al_{45}+(Q-h)\be_{45}
$$
or
$$
h = (\al+Q\be_{45})/\be.
$$
Therefore, the condition $Q>h>0$ is equivalent to
$$
Q > (\al+Q\be_{45})/\be > 0
$$
or
$$
Q>\max\{{\al/\be_{12}}; {-\al/\be_{45}}\}.
$$
Finally,
$$
T_{Eq} = T_1 = \al_{12}+(\al+Q\be_{45})\be_{12}/\be.
$$
\end{pf}
\qqed

The equilibrium in $N^+$ is described by seven cases in Lemma
\ref{L2}. It may be pointed out that these cases correspond to the
following situations in $N^+$:
 (a) the only path with non-vanishing flow is $P_3$;
 (b) the only path with non-vanishing flow is $P_2$;
 (c) the only path with non-vanishing flow is $P_1$;
 (d) the only path with vanishing flow is $P_1$;
 (e) the only path with vanishing flow is $P_2$;
 (f) the only path with vanishing flow is $P_3$;
 (g) no path has a vanishing flow.


Also, it is not difficult to see that some of the cases in Lemma
\ref{L2} are mutually exclusive, so the equilibrium in a
particular network $N^+$ is described by some of the the presented
seven cases. For example, if $\al_i=\be_i=1$ for $1\le i\le5$,
then the equilibrium is given by just one case (f).

\begin{lem}
\label{L2}
In the network $N^+$, the travel time at equilibrium is
as follows:

\begin{itemize}

\item [(a)] $T^+_{Eq} = \al_{135}+Q\be_{135}\;\;$ if $\;\; 0<Q\le
\min\{{\hat{\al}/\be_{35}}; {\bar{\al}/\be_{13}}\}$;

\item [(b)] $T^+_{Eq} = \al_{45}+Q\be_{45}\;\;$ if $\;\; 0<Q\le
\min\{{-\al/\be_{45}}; {-\bar{\al}/\be_{4}}\}$;

\item [(c)] $T^+_{Eq} = \al_{12}+Q\be_{12}\;\;$ if $\;\; 0<Q\le
\min\{{\al/\be_{12}}; {-\hat{\al}/\be_{2}}\}$;

\item [(d)] $T^+_{Eq} =
\al_{45}+Q\be_{45}-(\bar\al+Q\be_4)\be_4/\be_{134}\;\;$ if $\;\;
\max\{{\bar\al/\be_{13}}; {-\bar{\al}/\be_{4}}\} <Q\le \mu_1$;

\item [(e)] $T^+_{Eq} =
\al_{12}+Q\be_{12}-(\hat\al+Q\be_2)\be_2/\be_{235}\;\;$ if $\;\;
\max\{{\hat\al/\be_{35}}; {-\hat{\al}/\be_{2}}\} <Q\le \mu_2$;

\item [(f)] $T^+_{Eq} = \al_{12}+(\al+Q\be_{45})\be_{12}/\be\;\;$
if $\;\; Q > \max\{{\al/\be_{12}}; {-\al/\be_{45}}\}\;\; $ and
$$
{\cal B}_1 \ge {\hat\al\be_{14} + \bar\al\be_{25} \over Q};
$$

\item [(g)] $T^+_{Eq} = \al_{12}+(\al+Q\be_{45})\be_{12}/\be +
g{\cal B}_1/\be,\;\;$ where
$$
g = {\bar\al\be-\al\be_{14}-Q{\cal B}_1 \over
\be_3\be+\be_{14}\be_{25}},
$$
if $\;\; Q > \max\{\mu_1;\mu_2\}\;\;$  and
$$
{\cal B}_1 < {\hat\al\be_{14} + \bar\al\be_{25} \over Q}.
$$
\end{itemize}
\end{lem}

\begin{pf}
Let us denote $f_{ab}$ by $f$, and $f_{bc}$ by $g$. Then
$f_{ac}=Q-f$ and, using the conservation-of-flow constraints,
$f_{bd}=f-g$ and $f_{cd}=Q-f+g$ (see Figure 3). We have
$$
T_1= \al_{12}+f\be_{12} -g\be_2, \quad  T_2=
\al_{45}+Q\be_{45}-f\be_{45}+g\be_5 \quad\mbox{and}\quad T_3=
\al_{135}+Q\be_{5}+f(\be_{1}-\be_5)+g\be_{35},
$$
where $T_i$ is the travel time on the path $P_i$.

\unitlength=1.00mm
\begin{picture}(150,40)
\put(50,20){\circle*{2}} \put(70,30){\circle*{2}}
\put(70,10){\circle*{2}} \put(90,20){\circle*{2}}

\put(45,19){\makebox{$a$}} \put(93,19){\makebox{$d$}}
\put(69,32){\makebox{$b$}}\put(69,6){\makebox{$c$}}

 \put(56,27){\makebox{\small $f$}}
 \put(79,27){\makebox{\small $f-g$}}
 \put(52,11){\makebox{\small $Q-f$}}
 \put(78,11){\makebox{\small $Q-f+g$}}
 \put(72,19){\makebox{\small $g$}}


\qbezier (50,20)(70,30)(70,30) \qbezier
(50,20)(70,10)(70,10)\qbezier (70,30)(70,10)(70,10)\qbezier
(70,30)(90,20)(90,20)\qbezier (70,10)(90,20)(90,20)


\qbezier (68,28)(69,29.5)(69,29.5)
\qbezier(67.3,29.5)(69,29.5)(69,29.5)

\qbezier (67.3,10.5)(69,10.5)(69,10.5)
\qbezier(68,12)(69,10.5)(69,10.5)

\qbezier (88,18)(89,19.5)(89,19.5)
\qbezier(87.3,19.5)(89,19.5)(89,19.5)

\qbezier (87.3,20.5)(89,20.5)(89,20.5)
\qbezier(88,22)(89,20.5)(89,20.5)

\qbezier(70.8,12.5)(70,11)(70,11)
\qbezier(69.2,12.5)(70,11)(70,11)


\end{picture}

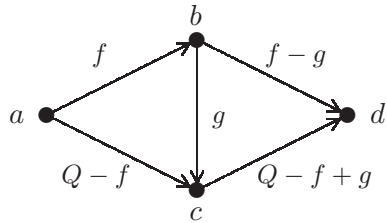
\vspace*{-.3cm} \centerline{{\bf Figure 3.} The link flows in
$N^+$.}

\bigskip

\2 \noindent {\bf\em Case (a):} The only path with non-vanishing
flow is $P_3$, i.e. $P_1$ and $P_2$ have a vanishing flow.
Therefore, $Q=f=g>0$ and
$$
T_1= \al_{12}+Q\be_{1}, \quad T_2= \al_{45}+Q\be_{5}
\quad\mbox{and}\quad T_3= \al_{135}+Q\be_{135}.
$$
At equilibrium, $T_1\ge T_3$ and $T_2\ge T_3$, i.e.
$$
\left\{\begin{array}{l}
    \hbox{$\al_{12}+Q\be_{1} \ge \al_{135}+Q\be_{135}$;} \\
    \hbox{$\al_{45}+Q\be_{5} \ge \al_{135}+Q\be_{135}$.} \\
\end{array}%
\right.
$$
Thus,
$$
0<Q\le \min\{{\hat{\al}/\be_{35}}; {\bar{\al}/\be_{13}}\}
$$
and
$$
T^+_{Eq} = T_3 = \al_{135}+Q\be_{135}.
$$

\2 \noindent {\bf\em Case (b):} The only path with non-vanishing
flow is $P_2$. Hence $Q>f=g=0$ and
$$
T_1= \al_{12}, \quad T_2= \al_{45}+Q\be_{45} \quad\mbox{and}\quad
T_3= \al_{135}+Q\be_{5}.
$$
At equilibrium, $T_1\ge T_2$ and $T_3\ge T_2$, i.e.
$$
\left\{\begin{array}{l}
    \hbox{$\al_{12} \ge \al_{45}+Q\be_{45}$;} \\
    \hbox{$\al_{135}+Q\be_{5} \ge \al_{45}+Q\be_{45}$.} \\
\end{array}%
\right.
$$
Thus,
$$
0<Q\le \min\{{-{\al}/\be_{45}}; {-\bar{\al}/\be_{4}}\}
$$
and
$$
T^+_{Eq} = T_2 = \al_{45}+Q\be_{45}.
$$

\2 \noindent {\bf\em Case (c):} The only path with non-vanishing
flow is $P_1$. Hence $Q=f>g=0$ and
$$
T_1= \al_{12}+Q\be_{12}, \quad T_2= \al_{45} \quad\mbox{and}\quad
T_3= \al_{135}+Q\be_{1}.
$$
At equilibrium, $T_2\ge T_1$ and $T_3\ge T_1$, i.e.
$$
\left\{\begin{array}{l}
    \hbox{$\al_{45} \ge \al_{12}+Q\be_{12}$;} \\
    \hbox{$\al_{135}+Q\be_{1} \ge \al_{12}+Q\be_{12}$.} \\
\end{array}%
\right.
$$
Thus,
$$
0<Q\le \min\{{{\al}/\be_{12}}; {-\hat{\al}/\be_{2}}\}
$$
and
$$
T^+_{Eq} = T_1 = \al_{12}+Q\be_{12}.
$$

\2 \noindent {\bf\em Case (d):} The only path with vanishing flow
is $P_1$. We obtain $Q>f=g>0$ and
$$
T_1= \al_{12}+f\be_{1}, \quad T_2= \al_{45}+Q\be_{45}-f\be_{4}
\quad\mbox{and}\quad T_3= \al_{135}+Q\be_{5}+f\be_{13}.
$$
At equilibrium, $T_2= T_3$ and $T_1\ge T_2$, i.e.
$$
\left\{\begin{array}{l}
    \hbox{$\al_{45}+Q\be_{45}-f\be_{4} = \al_{135}+Q\be_{5}+f\be_{13}$;}\\
    \hbox{$\al_{12}+f\be_{1} \ge \al_{45}+Q\be_{45}-f\be_{4}$;}\\
\end{array}%
\right. \quad\mbox{or}\quad
 \left\{\begin{array}{l}
    \hbox{$f = (\bar\al+Q\be_{4})/\be_{134};$} \\
    \hbox{$f \ge (\al+Q\be_{45})/\be_{14}$.} \\
\end{array}%
\right.
$$
Therefore,
$$
(\bar\al+Q\be_{4})/\be_{134}  \ge (\al+Q\be_{45})/\be_{14}
$$
or
$$
Q(\be_{45}\be_{134}-\be_{4}\be_{14}) \le \bar\al\be_{14} -
\al\be_{134}.
$$
Note that $\be_{45}\be_{134}-\be_{4}\be_{14} =
\be_{3}\be_{45}+\be_{5}\be_{14}$ and, using $\bar\al=\al+\hat\al$,
 \beq \label{eq1}
\bar\al\be_{14} - \al\be_{134} = \hat\al\be_{14}-\al\be_3.
 \eeq
Hence
$$
Q \le {\hat{\al}\be_{14} - \al\be_{3} \over
\be_3\be_{45}+\be_5\be_{14}},
$$
i.e. $Q\le \mu_1$. Now, taking into account that $Q>f>0$, the
following is obtained:
$$
Q > (\bar\al+Q\be_{4})/\be_{134}  > 0,
$$
which is equivalent to
$$
Q > \max\{{\bar\al/\be_{13}}; {-\bar{\al}/\be_{4}}\}.
$$
Finally,
$$
T^+_{Eq} = T_2 =
\al_{45}+Q\be_{45}-(\bar\al+Q\be_4)\be_4/\be_{134}.
$$

\2 \noindent {\bf\em Case (e):} The only path with vanishing flow
is $P_2$. We obtain $Q=f>g>0$ and
$$
T_1= \al_{12}+Q\be_{12}-g\be_2, \quad T_2= \al_{45}+g\be_{5}
\quad\mbox{and}\quad T_3= \al_{135}+Q\be_{1}+g\be_{35}.
$$
At equilibrium, $T_1= T_3$ and $T_2\ge T_3$, i.e.
$$
\left\{\begin{array}{l}
    \hbox{$\al_{12}+Q\be_{12}-g\be_2= \al_{135}+Q\be_{1}+g\be_{35}$;}\\
    \hbox{$\al_{45}+g\be_{5} \ge \al_{135}+Q\be_{1}+g\be_{35}$;}\\
\end{array}%
\right. \quad\mbox{or}\quad
 \left\{\begin{array}{l}
    \hbox{$g = (\hat\al+Q\be_{2})/\be_{235};$} \\
    \hbox{$g \le (\bar\al-Q\be_{1})/\be_{3}$.} \\
\end{array}%
\right.
$$
Therefore,
$$
(\hat\al+Q\be_{2})/\be_{235}  \le (\bar\al-Q\be_{1})/\be_{3}
$$
or
$$
Q(\be_{2}\be_{3}+\be_{1}\be_{235}) \le \bar\al\be_{235} -
\hat\al\be_{3}.
$$
Now, rearranging the left-hand side, and using
$\hat\al=\bar\al-\al$, we obtain
$$
Q(\be_{1}\be_{25}+\be_{3}\be_{12}) \le \bar\al\be_{25}
+\al\be_{3}.
$$
Thus, $Q \le \mu_2$. Taking into account that $Q>g>0$, we obtain
$$
Q > (\hat\al+Q\be_{2})/\be_{235}  > 0,
$$
which is equivalent to
$$
Q > \max\{{\hat\al/\be_{35}}; {-\hat{\al}/\be_{2}}\}.
$$
Finally,
$$
T^+_{Eq} = T_1 =
\al_{12}+Q\be_{12}-(\hat\al+Q\be_2)\be_2/\be_{235}.
$$

\2 \noindent {\bf\em Case (f):} The only path with vanishing flow
is $P_3$. We obtain $Q>f>g=0$ and
$$
T_1= \al_{12}+f\be_{12}, \quad T_2= \al_{45}+Q\be_{45}-f\be_{45}
\quad\mbox{and}\quad T_3= \al_{135}+Q\be_{5}+f(\be_{1}-\be_5).
$$
At equilibrium, $T_1= T_2$ and $T_3\ge T_2$, i.e.
$$
\left\{\begin{array}{l}
    \hbox{$\al_{12}+f\be_{12} = \al_{45}+Q\be_{45}-f\be_{45}$;}\\
    \hbox{$\al_{135}+Q\be_{5}+f(\be_{1}-\be_5) \ge  \al_{45}+Q\be_{45}-f\be_{45}$;}\\
\end{array}%
\right. \quad\mbox{or}\quad
 \left\{\begin{array}{l}
    \hbox{$f = (\al+Q\be_{45})/\be;$} \\
    \hbox{$f \ge (\bar\al+Q\be_{4})/\be_{14}$.} \\
\end{array}%
\right.
$$
Therefore,
$$
 (\al+Q\be_{45})/\be  \ge (\bar\al+Q\be_{4})/\be_{14}
$$
or
$$
Q(\be_{14}\be_{45}-\be\be_{4}) \ge \bar\al\be - \al\be_{14}.
$$
It is not difficult to see that $Q(\be_{14}\be_{45}-\be\be_{4})
=Q{\cal B}_1$ and, using $\al=\bar\al-\hat\al$,
 \beq \label{eq2}
 \bar\al\be -
\al\be_{14} =
\bar\al(\be_{14}+\be_{25})-(\bar\al-\hat\al)\be_{14}=
\hat\al\be_{14}+\bar\al\be_{25}.
 \eeq
 Thus,
$$
{\cal B}_1 \ge {\hat\al\be_{14}+\bar\al\be_{25}\over Q}.
$$
Taking into account that $Q>f>0$, we obtain
$$
Q > (\al+Q\be_{45})/\be > 0,
$$
which is equivalent to
$$
Q > \max\{{\al/\be_{12}}; {-\al/\be_{45}}\}.
$$
Finally,
$$
T^+_{Eq} = T_1 = \al_{12}+(\al+Q\be_{45})\be_{12}/\be.
$$

\2 \noindent {\bf\em Case (g):} No path has a vanishing flow, and
so $Q>f>g>0$. At equilibrium, $T_1= T_2$ and $T_2 = T_3$, i.e.
$$
\left\{\begin{array}{l}
    \hbox{$f\be = \al+g\be_{25}+Q\be_{45}$;}\\
    \hbox{$f\be_{14}+g\be_{3} =  \bar\al+Q\be_{4}$;}\\
\end{array}%
\right. \quad\mbox{or}\quad
 \left\{\begin{array}{l}
    \hbox{$f = (\al+g\be_{25}+Q\be_{45})/\be;$} \\
    \hbox{$g = \left(\bar\al\be-\al\be_{14}-Q(\be_{45}\be_{14}-\be_4\be)\right) / (\be_3\be+\be_{14}\be_{25})$.} \\
\end{array}%
\right.
$$
It is easy to see that $\be_{45}\be_{14}-\be_4\be = {\cal B}_1$,
so
$$
g = {\bar\al\be-\al\be_{14}-Q{\cal B}_1 \over
\be_3\be+\be_{14}\be_{25}}.
$$

The condition $g>0$ is equivalent to
$\bar\al\be-\al\be_{14}-Q{\cal B}_1 > 0 $ or $ Q{\cal B}_1 <
\bar\al\be-\al\be_{14}$. Using (\ref{eq2}), we obtain
$$
{\cal B}_1 < {\hat\al\be_{14} + \bar\al\be_{25} \over Q}.
$$

The condition $f>g$ can be written as
$(\al+g\be_{25}+Q\be_{45})/\be > g$ or $g <
(\al+Q\be_{45})/\be_{14}$. Hence
$$
{\bar\al\be-\al\be_{14}-Q{\cal B}_1 \over
\be_3\be+\be_{14}\be_{25}} < {\al+Q\be_{45} \over \be_{14}},
$$
which is equivalent to
$$
Q({\cal B}_1 \be_{14}  + \be_{45}\be_{25}\be_{14} +
\be_{3}\be_{45}\be) > \bar\al\be_{14}\be - \al(\be^2_{14} +
\be_{3}\be + \be_{25}\be_{14}).
$$
It is easy to see that ${\cal B}_1 \be_{14} +
\be_{45}\be_{25}\be_{14} = \be_5\be_{14}\be$ and $\be^2_{14} +
\be_{3}\be + \be_{25}\be_{14} = \be_{134}\be$. Therefore,
$$
Q(\be_{3}\be_{45}\be + \be_{5}\be_{14}\be)
> \bar\al\be_{14}\be - \al\be_{134}\be \quad \mbox{or}\quad Q(\be_{3}\be_{45} + \be_{5}\be_{14})
> \bar\al\be_{14} - \al\be_{134}.
$$
Using (\ref{eq1}), we obtain $Q>\mu_1$.

Now let us consider the condition $Q>f$, which can be written as
$Q > (\al+g\be_{25}+Q\be_{45})/\be$ or $(Q\be_{12}-\al)/\be_{25}
> g$. Hence
$$
{Q\be_{12}-\al\over \be_{25}} > {\bar\al\be-\al\be_{14}-Q{\cal
B}_1 \over \be_3\be+\be_{14}\be_{25}}
$$
or
$$
Q(\be_{12}\be_{25}\be_{14}+\be_{25}{\cal B}_1 +
\be_{3}\be_{12}\be) > \bar\al\be_{25}\be+\al\be_{3}\be.
$$
It is not difficult to check that
$\be_{12}\be_{25}\be_{14}+\be_{25}{\cal B}_1 = \be_{1}\be_{25}\be$
and hence
$$
Q( \be_{1}\be_{25}\be + \be_{3}\be_{12}\be)
> \bar\al\be_{25}\be+\al\be_{3}\be,
$$
i.e. $Q>\mu_2$.

Finally,
$$
T^+_{Eq} = T_1 = \al_{12}+f\be_{12} -g\be_2 =
\al_{12}+(\al+Q\be_{45})\be_{12}/\be +
g(\be_{25}\be_{12}/\be-\be_2).
$$
It is easy to see that $\be_{25}\be_{12}/\be-\be_2 = {\cal
B}_1/\be$, and so
$$
T^+_{Eq} = \al_{12}+(\al+Q\be_{45})\be_{12}/\be + g{\cal B}_1/\be.
$$
\end{pf}
\qqed


\section{Braess' Paradox in $N$ and $N^+$}

Braess' paradox is said to occur in the network configuration
$N/N^+$ for a given total flow $Q$ if
$$
T^+_{Eq} > T_{Eq},
$$
where  $T_{Eq}$ and $T^+_{Eq}$ are travel times at equilibria in
$N$ and $N^+$, respectively. Thus, Braess' paradox describes a
situation when adding a new link to a network makes a general
performance worse.

Our first result describes all possible situations when Braess'
paradox may occur in $N/N^+$ in terms of their paths. In fact, it
says that Braess' paradox may occur in $N/N^+$ only if, at
equilibria, both $P_1$ and $P_2$ have a non-vanishing flow in $N$,
and $P_3$ has a non-vanishing flow in $N^+$.

\2

\noindent{\bf Mega-Theorem.} \label{MT} {\it Braess' paradox may
occur in $N/N^+$ in the following cases only:
\begin{itemize}
    \item [(a)] At equilibria, both $N$ and $N^+$ have no paths with
    vanishing flow.
    \item [(b)] At equilibria, $N$ has no path with vanishing flow, and  $P_3$ is the only
    path with non-vanishing flow in  $N^+$.
    \item [(c)] At equilibria, $N$ has no path with vanishing flow, and  $P_1$ is the only
    path with vanishing flow in  $N^+$.
    \item [(d)] At equilibria, $N$ has no path with vanishing flow, and  $P_2$ is the only
    path with vanishing flow in  $N^+$.
\end{itemize}}
\2

In the following theorems, we formulate the necessary and
sufficient conditions for the existence of the paradox for all the
cases of Mega-Theorem. These theorems are proved in Section
\ref{Proof}.

\begin{thm}
\label{T1} Suppose that at equilibria both $N$ and $N^+$ have no
paths with vanishing flow. Then Braess' paradox occurs in $N/N^+$
if and only if the Braess number ${\cal B}_1$ is positive and
$$
\max\left\{{\al\over\be_{12}}; {-\al\over\be_{45}}; \mu_1;
\mu_2\right\} < Q < {\hat\al\be_{14} + \bar\al\be_{25} \over {\cal
B}_1}.
$$
\end{thm}

\begin{thm}
\label{T2} Suppose that  at equilibria $N$ has no path with
vanishing flow and $P_3$ is the only  path with non-vanishing flow
in $N^+$. Then Braess' paradox occurs in $N/N^+$ if and only if
the Braess number ${\cal B}_2$ is positive and
$$
\max\left\{{\al\over\be_{12}}; {-\al\over\be_{45}};{\hat\al
\be_{45} + \bar\al \be_{12}\over {\cal B}_2 }\right\} < Q \le \min
\left\{{\hat\al\over\be_{35}}; {\bar\al\over\be_{13}}\right\}.
$$
\end{thm}

\begin{thm}
\label{T3} Suppose that  at equilibria $N$ has no path with
vanishing flow and $P_1$ is the only    path with vanishing flow
in  $N^+$. Then Braess' paradox occurs in $N/N^+$ if and only if
the Braess number ${\cal B}_3$ is positive and
$$
\max\left\{{\al\over\be_{12}}; {-\al\over\be_{45}}; {\bar\al\over
\be_{13}}; {-\bar\al\over\be_4};  { \bar\al\be_4\be -
\al\be_{134}\be_{45} \over  {\cal B}_3 } \right\} < Q \le \mu_1.
$$
\end{thm}

\begin{thm}
\label{T4} Suppose that  at equilibria $N$ has no path with
vanishing flow and $P_2$ is the only    path with vanishing flow
in  $N^+$. Then Braess' paradox occurs in $N/N^+$ if and only if
the Braess number ${\cal B}_4$ is positive and
$$
\max\left\{{\al\over\be_{12}}; {-\al\over\be_{45}}; {\hat\al\over
\be_{35}}; {-\hat\al\over\be_2}; { \hat\al\be_2\be +
\al\be_{235}\be_{12} \over  {\cal B}_4 } \right\} < Q \le \mu_2.
$$
\end{thm}

It might be pointed out that if ${\cal B}_1 \ge 0$, then ${\cal
B}_2$, ${\cal B}_3$ and ${\cal B}_4$ are positive numbers, because
$$
{\cal B}_2 = \be_{12}\be_{13} + \be_{35}\be_{45} + {\cal B}_1,
$$
$$
{\cal B}_3 = \be^2_{5}\be_{134} + \be_4(\be_{3}\be_{455} +
\be_{5}\be_{14} + {\cal B}_1),
$$
$$
{\cal B}_4 = \be^2_{1}\be_{235} + \be_2(\be_{1}\be_{335} +
\be_{2}\be_{13} + {\cal B}_1).
$$
Moreover, Theorems \ref{T3} and \ref{T4} are mutually exclusive in
the sense that they cannot provide intervals for $Q$
simultaneously. This is true because the inequalities
$\bar\al/\be_{13} < \mu_1$ and $\hat\al/\be_{35} < \mu_2$ are
inconsistent. Note also that if, for example, Theorems
\ref{T1}--\ref{T3} provide non-empty intervals for $Q$, then the
interval with highest values of $Q$ is given by Theorem \ref{T1},
the interval with smallest values of $Q$ is provided by Theorem
\ref{T2}, and Theorem \ref{T3} yields the interval with mid-range
values of $Q$.

Note that the original assumption $\be_i>0$ for all $i$ can be
relaxed by allowing $\be_i=0$ for some $i$. This can be done by
introducing $+\infty$ and $-\infty$ when a non-zero number is
divided by zero. For example, let us consider Arnott--Small's
example \cite{Arn}:
$$
\al_1=\al_5=0, \;\; \al_2=\al_4=15, \;\; \al_3=7.5, \;\;
\be_1=\be_5=0.01, \;\; \be_2=\be_3=\be_4=0.
$$
We obtain
$$
\al=0, \;\; \bar\al=\hat\al=7.5, \;\; \be=0.02, \;\; {\cal
B}_1=10^{-4}, \;\; \mu_1=\mu_2=750.
$$
By Theorems \ref{T1} and \ref{T2}, Braess' paradox occurs if
$$
500 < Q < 1500,
$$
while Theorems \ref{T3} and \ref{T4} provide no intervals for $Q$.
In calculating the lower bounds for $Q$ in these theorems we have
division by zero, but this problem is overcome by putting
$-7.5/0=-\infty$.

\section{Symmetrical/Asymmetrical Networks and the Pseudo-Paradox}

Let us consider the classical case of a symmetrical network
presented by Braess \cite{Bra} and discussed in Pas and Principio
\cite{Pas} and other papers. Using our notation, it is a
particular case of the network configuration $N/N^+$ when time
functions are symmetrical for links that do not share nodes with
each other ($(a,b)$ and $(c,d)$, $(a,c)$ and $(b,d)$), the free
flow travel times for the links $(a,b)$ and $(c,d)$ are equal to
zero, and the delay parameter for $(b,c)$ is equal to the delay
parameter of the links $(b,d)$ and $(a,c)$, i.e.
$$
\al_1=\al_5=0, \;\; \al_2=\al_4, \;\; \be_1=\be_5, \;\;
\be_2=\be_3=\be_4.
$$
Also, it is assumed that $\al_2>\al_3$ and $\be_1>\be_2$. This
network configuration is denoted by $M/M^+$ (see Figure 4).

\unitlength=1.00mm
\begin{picture}(150,40)
\put(50,20){\circle*{2}} \put(70,30){\circle*{2}}
\put(70,10){\circle*{2}} \put(90,20){\circle*{2}}

\put(45,19){\makebox{$a$}} \put(93,19){\makebox{$d$}}
\put(69,32){\makebox{$b$}}\put(69,6){\makebox{$c$}}

\put(54,27){\makebox{\small $0,\be_1$}}
\put(77,27){\makebox{\small $\al_2,\be_2$}}
\put(54,10.5){\makebox{\small $\al_2,\be_2$}}
\put(77,10.5){\makebox{\small $0,\be_1$}}
\put(71,19){\makebox{\small $\al_3,\be_2$}}


\qbezier (50,20)(70,30)(70,30) \qbezier
(50,20)(70,10)(70,10)\qbezier (70,30)(70,10)(70,10)\qbezier
(70,30)(90,20)(90,20)\qbezier (70,10)(90,20)(90,20)


\qbezier (68,28)(69,29.5)(69,29.5)
\qbezier(67.3,29.5)(69,29.5)(69,29.5)

\qbezier (67.3,10.5)(69,10.5)(69,10.5)
\qbezier(68,12)(69,10.5)(69,10.5)

\qbezier (88,18)(89,19.5)(89,19.5)
\qbezier(87.3,19.5)(89,19.5)(89,19.5)

\qbezier (87.3,20.5)(89,20.5)(89,20.5)
\qbezier(88,22)(89,20.5)(89,20.5)

\qbezier(70.8,12.5)(70,11)(70,11)
\qbezier(69.2,12.5)(70,11)(70,11)


\end{picture}

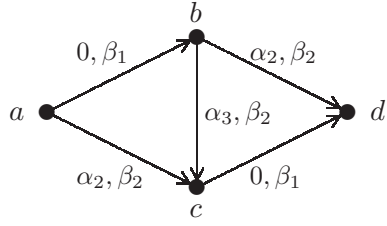
\vspace*{-.3cm} \centerline{{\bf Figure 4.} The symmetric network
$M^+$ ($\al_2>\al_3$ and $\be_1>\be_2$).}

\bigskip

Pas and Principio \cite{Pas} determined the occurrence of Braess'
paradox in the symmetrical network configuration $M/M^+$. This
result follows directly from Theorems \ref{T1}--\ref{T4}:

\begin{cor}
[\cite{Pas}] \label{PP}
Braess' paradox occurs in $M/M^+$ if and
only if
$$
{2(\al_2-\al_3)\over 3\be_1+\be_2} < Q < {2(\al_2-\al_3)\over
\be_1-\be_2}.
$$
\end{cor}

\begin{pf}
For the network configuration $M/M^+$, we have $$\al=0,\;\;
\bar\al=\hat\al=\al_2-\al_3,\;\;
\mu_1=\mu_2=(\al_2-\al_3)/\be_{12}, \;\; {\cal
B}_1=\be_1^2-\be_2^2.
$$
The Braess number ${\cal B}_1$ is positive because $\be_1>\be_2$
in $M/M^+$. Under the conditions of Theorem \ref{T1}, Braess'
paradox occurs if
$$
{\al_2-\al_3\over\be_{12}} < Q <
{2\be_{12}(\al_2-\al_3)\over\be_1^2-\be_2^2} =
{2(\al_2-\al_3)\over\be_1-\be_2}.
$$
Now,
$$
{\cal B}_2= \be_1^2+2\be^2_{12}-\be_2^2=\be_{12}(3\be_1+\be_2),
$$
which is a positive number. Therefore, by Theorem \ref{T2},
$$
{2(\al_2-\al_3)\over 3\be_1+\be_2} < Q \le
{\al_2-\al_3\over\be_{12}}.
$$
Note that the lower bound is less than the upper bound because
$\be_1>\be_2$. Thus, the above inequalities can be written
together as
$$
{2(\al_2-\al_3)\over 3\be_1+\be_2} < Q <
{2(\al_2-\al_3)\over\be_1-\be_2}.
$$
The upper and lower bounds of Theorems \ref{T3} and \ref{T4}
provide no intervals for $Q$.
\end{pf}
\qqed

Now let $\tilde{M}/\tilde{M}^+$ denote the above network
configuration $M/M^+$ without the assumption that $\al_2>\al_3$
and $\be_1>\be_2$. A proof similar to that of Corollary \ref{PP}
shows that Braess' paradox occurs in $\tilde{M}/\tilde{M}^+$ in
the following cases only:

(a) If $\be_1>\be_2$ and
$$
{\al_2-\al_3\over \be_{12}} < Q < {2(\al_2-\al_3)\over
\be_1-\be_2};
$$

(b) If
$$
{2(\al_2-\al_3)\over 3\be_1+\be_2} < Q \le {\al_2-\al_3\over
\be_{12}}.
$$

The both cases imply that $\al_2>\al_3$. Another implicit
relationship is obtained from (b) if we require that
$$
{2(\al_2-\al_3)\over 3\be_1+\be_2} <  {\al_2-\al_3\over \be_{12}},
$$
which is equivalent to $\be_1>\be_2$. Thus, even though the
network configuration $\tilde{M}/\tilde{M}^+$ extends $M/M^+$, the
conditions for the occurrence of Braess' paradox are the same.

Let us further extend the network configuration
$\tilde{M}/\tilde{M}^+$ by allowing any non-negative free flow
travel time $\al_1$ for the links $(a,b)$ and $(c,d)$ and any
positive delay parameter $\be_3$ for the link $(b,c)$. In other
words, the symmetrical network configuration $S/S^+$ of Figure 5
is obtained from $N/N^+$ when time functions are symmetrical for
links that do not share nodes with each other:
$$
\al_1=\al_5, \;\; \al_2=\al_4, \;\; \be_1=\be_5, \;\; \be_2=\be_4.
$$

\unitlength=1.00mm
\begin{picture}(150,40)
\put(50,20){\circle*{2}} \put(70,30){\circle*{2}}
\put(70,10){\circle*{2}} \put(90,20){\circle*{2}}

\put(45,19){\makebox{$a$}} \put(93,19){\makebox{$d$}}
\put(69,32){\makebox{$b$}}\put(69,6){\makebox{$c$}}

\put(54,27){\makebox{\small $\al_1,\be_1$}}
\put(77,27){\makebox{\small $\al_2,\be_2$}}
\put(54,11){\makebox{\small $\al_2,\be_2$}}
\put(77,11){\makebox{\small $\al_1,\be_1$}}
\put(71,19){\makebox{\small $\al_3,\be_3$}}


\qbezier (50,20)(70,30)(70,30) \qbezier
(50,20)(70,10)(70,10)\qbezier (70,30)(70,10)(70,10)\qbezier
(70,30)(90,20)(90,20)\qbezier (70,10)(90,20)(90,20)


\qbezier (68,28)(69,29.5)(69,29.5)
\qbezier(67.3,29.5)(69,29.5)(69,29.5)

\qbezier (67.3,10.5)(69,10.5)(69,10.5)
\qbezier(68,12)(69,10.5)(69,10.5)

\qbezier (88,18)(89,19.5)(89,19.5)
\qbezier(87.3,19.5)(89,19.5)(89,19.5)

\qbezier (87.3,20.5)(89,20.5)(89,20.5)
\qbezier(88,22)(89,20.5)(89,20.5)

\qbezier(70.8,12.5)(70,11)(70,11)
\qbezier(69.2,12.5)(70,11)(70,11)


\end{picture}

\vspace*{-.3cm} \centerline{{\bf Figure 5.} The symmetric network
$S^+$.}

\bigskip

\begin{cor}
\label{C2} Braess' paradox occurs in the symmetrical network
configuration $S/S^+$ if and only if $\be_1>\be_2$ and
 \beq \label{eq0}
 {2(\al_2-\al_{13})\over 3\be_{1}+2\be_3-\be_2} < Q <
{2(\al_2-\al_{13})\over \be_1-\be_2}.
 \eeq
\end{cor}

\begin{pf}
For the network configuration $S/S^+$, we have $$\al=0,\;\;
\bar\al=\hat\al=\al_2-\al_{13},\;\;
\mu_1=\mu_2=(\al_2-\al_{13})/\be_{13}, \;\; {\cal
B}_1=\be_1^2-\be_2^2 = \be_{12}(\be_1-\be_2).
$$
By Theorem \ref{T1}, Braess' paradox occurs if ${\cal B}_1$ is
positive, i.e.  $\be_1>\be_2$, and
$$
{\al_2-\al_{13}\over\be_{13}} < Q <
{2(\al_2-\al_{13})\over\be_1-\be_2} .
$$
Now,
$$
{\cal B}_2=
\be_{131}2\be_{12}-\be_{12}^2=\be_{12}(3\be_1+2\be_3-\be_2).
$$
Therefore, by Theorem \ref{T2}, Braess' paradox occurs if
$3\be_1+2\be_3>\be_2$ and
$$
{2(\al_2-\al_{13})\over 3\be_{1}+2\be_3-\be_2} < Q \le
{\al_2-\al_{13}\over \be_{13}}.
$$
This implies that $\al_2>\al_{13}$. Also, it is easy to see that
the lower bound is less than the upper bound only if
$\be_1>\be_2$, which is stronger than $3\be_1+2\be_3>\be_2$. Thus,
the above inequalities can be written together as
$$
{2(\al_2-\al_{13})\over 3\be_{1}+2\be_3-\be_2} < Q <
{2(\al_2-\al_{13})\over \be_1-\be_2}.
$$
The upper and lower bounds of Theorems \ref{T3} and \ref{T4}
provide no intervals for $Q$.
\end{pf}
\qqed

In Corollary \ref{C2} there is an implicit assumption that
$\al_2>\al_{13}$ because $Q$ is a positive number (if $\al_2\le
\al_{13}$, then (\ref{eq0}) provides no interval for $Q$). We will
see in Corollary \ref{C6} what is happening with the times at
equilibria in $S/S^+$ if $Q$ exceeds the upper bound in
(\ref{eq0}), where $\al_2>\al_{13}$ and $\be_1>\be_2$.

\unitlength=1.00mm
\begin{picture}(150,40)
\put(50,20){\circle*{2}} \put(70,30){\circle*{2}}
\put(70,10){\circle*{2}} \put(90,20){\circle*{2}}

\put(45,19){\makebox{$a$}} \put(93,19){\makebox{$d$}}
\put(69,32){\makebox{$b$}}\put(69,6){\makebox{$c$}}

\put(54,27){\makebox{\small $\al_1,\be_1$}}
\put(77,27){\makebox{\small $\al_2,\be_2$}}
\put(54,11){\makebox{\small $\al_1,\be_1$}}
\put(77,11){\makebox{\small $\al_2,\be_2$}}
\put(71,19){\makebox{\small $\al_3,\be_3$}}


\qbezier (50,20)(70,30)(70,30) \qbezier
(50,20)(70,10)(70,10)\qbezier (70,30)(70,10)(70,10)\qbezier
(70,30)(90,20)(90,20)\qbezier (70,10)(90,20)(90,20)


\qbezier (68,28)(69,29.5)(69,29.5)
\qbezier(67.3,29.5)(69,29.5)(69,29.5)

\qbezier (67.3,10.5)(69,10.5)(69,10.5)
\qbezier(68,12)(69,10.5)(69,10.5)

\qbezier (88,18)(89,19.5)(89,19.5)
\qbezier(87.3,19.5)(89,19.5)(89,19.5)

\qbezier (87.3,20.5)(89,20.5)(89,20.5)
\qbezier(88,22)(89,20.5)(89,20.5)

\qbezier(70.8,12.5)(70,11)(70,11)
\qbezier(69.2,12.5)(70,11)(70,11)


\end{picture}

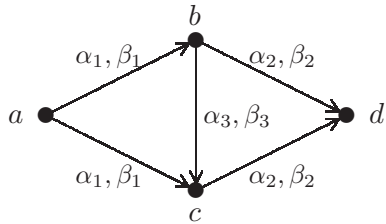
\vspace*{-.3cm} \centerline{{\bf Figure 6.} The asymmetric network
$A^+$.}

\bigskip

The asymmetrical network configuration $A/A^+$ of Figure 6 is
obtained from $N/N^+$ when time functions are same for links that
share the origin and the destination:
$$
\al_1=\al_4, \;\; \al_2=\al_5, \;\; \be_1=\be_4, \;\; \be_2=\be_5.
$$

\begin{cor}
\label{C3} Braess' paradox cannot occur in the asymmetrical
network configuration $A/A^+$.
\end{cor}

\begin{pf}
For the network configuration $A/A^+$, we have $$\al=0,\;\;
\bar\al=\hat\al=-\al_{3},\;\; {\cal B}_1=\be_1\be_2-\be_2\be_1 =0,
\;\; \mu_1\le 0, \;\; \mu_2\le 0.
$$
It is easy to see that Theorems \ref{T1}--\ref{T4} provide no
intervals, so Braess' paradox is impossible.
\end{pf}
\qqed

Using Lemmas \ref{L1} and  \ref{L2} we see that the equilibria in
$A/A^+$ are described by the cases (c) and (f), respectively, i.e.
in  $A$ no path has a vanishing flow, and in  $A^+$ the only path
with vanishing flow is $P_3$. Also, the flow $Q$ is distributed
evenly between $P_1$ and $P_2$ and
$$
T_{Eq}=T^+_{Eq} = \al_{12}+0.5\be_{12}Q.
$$
Thus, the travel times at equilibria in  $A$ and $A^+$ are equal
for any $Q$. This observation is important because adding a new
link to $A$ does not improve the general performance, even though
Braess' paradox is not occurring.

We say that the pseudo-paradox occurs in the network configuration
$N/N^+$ if
$$
T^+_{Eq} = T_{Eq}
$$
for an interval of values for $Q$ (as opposite to a single point).
In other words, we exclude single values of $Q$ when going from
the situation ``Braess' paradox does not occur'' to ``Braess'
paradox occurs''.

Thus, the pseudo-paradox describes a situation when adding a new
link to a network does not change the general performance for a
range of values of the total flow. As seen above, the
pseudo-paradox occurs in the asymmetrical network configuration
$A/A^+$ for any $Q>0$. In our view, the pseudo-paradox is a more
common phenomenon than Braess' paradox.

\begin{cor}
\label{C4} For the network configuration $N/N^+$, the
pseudo-paradox occurs if
\begin{itemize}
    \item [(a)] $0<Q\le \min\{-\al/\be_{45}; -\bar\al/\be_4\}$;
    \item [(b)] $0<Q\le \min\{\al/\be_{12}; -\hat\al/\be_2\}$;
    \item [(c)] $Q > \max\{\al/\be_{12}; -\al/\be_{45}\}\;\;$ and
    $\;\;
    Q{\cal B}_1 \ge \hat\al\be_{14} + \bar\al\be_{25};
    $
    \item [(d)] ${\cal B}_1=0$,
    $\;\;\hat\al\be_{14}+\bar\al\be_{25}>0\;\;$ and
    $\;\;
    Q > \max\{\al/\be_{12}; -\al/\be_{45}; \mu_1; \mu_2\}.
    $
\end{itemize}
\end{cor}

\begin{pf}
The first three cases follow immediately from the cases (a,b),
(b,c), (c,f) of the proof of theorems as discussed in the first
paragraph of the next section. The last case is similar to the
case (c,g), where $T^+_{Eq} = T_{Eq}$ must be satisfied. This is
equivalent to $g{\cal B}_1/\be=0$. Since $g>0$, we obtain ${\cal
B}_1=0$, and therefore $\hat\al\be_{14}+\bar\al\be_{25}>0$. In
addition, there are lower bounds for $Q$ in Lemma \ref{L1} (c) and
Lemma \ref{L2} (g).
\end{pf}
\qqed

The application of Corollary \ref{C4} (c) to the asymmetrical
network configuration $A/A^+$ confirms the above observation:

\begin{cor}
\label{C5} The pseudo-paradox occurs in the asymmetrical network
configuration $A/A^+$ for any $Q>0$.
\end{cor}

By Corollary \ref{C2}, Braess' paradox occurs in the symmetrical
network configuration if $\be_1>\be_2$, $\al_2 > \al_{13}$, and
the total flow $Q$ is between the lower and upper bounds in
(\ref{eq0}). Corollary \ref{C4} (c) allows us to see what is
happening with the times at equilibria if $Q$ exceeds the upper
bound:

\begin{cor}
\label{C6} Suppose that Braess' paradox occurs in the symmetrical
network configuration $S/S^+$, i.e. $\be_1>\be_2$ and $\al_2 >
\al_{13}$. Then $S/S^+$ is experiencing the pseudo-paradox for any
$$
Q \ge {2(\al_2-\al_{13})\over \be_1-\be_2}.
$$
\end{cor}

\begin{pf}
We know that $\al=0$, $\bar\al=\hat\al=\al_2-\al_{13}$ and ${\cal
B}_1=\be_{12}(\be_1-\be_2)>0$. Therefore, the second inequality in
Corollary \ref{C4} (c) is equivalent to
$$Q\be_{12}(\be_1-\be_2)
\ge 2(\al_2-\al_{13})\be_{12}, $$
as required.
\end{pf}
\qqed

Thus, under the conditions of Corollary \ref{C6}, some improvement
in $S/S^+$ is only possible if $Q$ is less than the lower bound in
(\ref{eq0}), followed by Braess' paradox until $Q$ reaches the
upper bound in (\ref{eq0}), followed by the pseudo-paradox for
larger values of $Q$.


We conclude this section with a numerical example. Let the network
configuration $N/N^+$ have the following parameters:
$$
\al_1=2, \;\; \al_2=36, \;\; \al_3=6, \;\; \al_4=40, \;\; \al_5=2,
$$
$$
\be_1=30, \;\; \be_2=32, \;\; \be_3=3, \;\; \be_4=8, \;\;
\be_5=19. \;\;
$$
We obtain
$$
\al=4, \;\; \bar\al=32,\;\; \hat\al=28,\;\; \be=89,\;\;
\mu_1=1.31,\;\; \mu_2=0.96,
$$
$$
{\cal B}_1=314, \;\; {\cal B}_2=2954, \;\;{\cal B}_3=24193,
\;\;{\cal B}_4=116440.
$$
Braess' paradox occurs in this network in the following cases:
$$
1.31<Q<8.59 \;\; \mbox{by Theorem \ref{T1}},
$$
$$ 0.93<Q\le 0.97 \;\; \mbox{by Theorem \ref{T2}},
$$
$$
0.97<Q\le 1.31 \;\; \mbox{by Theorem \ref{T3}}.
$$
Theorem \ref{T4} produces no interval. Therefore, Braess' paradox
occurs if and only if
$$
0.93<Q<8.59.
$$
Note that rounded numbers are used instead of exact values. By
Corollary \ref{C4} (c), the pseudo-paradox occurs if
$$
Q\ge 8.59.
$$
Thus, some improvement in the network when adding the link $(b,c)$
only occurs for small values of $Q$ ($Q<0.93$). The extent of this
improvement and Braess' paradox can be seen from the equilibrium
functions found by Lemmas \ref{L1} and \ref{L2}:
$$
T_{Eq} = \left\{\begin{array}{l}
    \hbox{$38+62Q\;\;$ if $\;\;0<Q\le0.06$;} \\
    \hbox{$40.79+18.81Q\;\;$ if $\;\;Q>0.06$;} \\
\end{array}%
\right.
$$
and
$$
T^+_{Eq} = \left\{\begin{array}{l}
    \hbox{$10+52Q\;\;$ if $\;\;0<Q\le 0.97$;} \\
    \hbox{$35.76+25.44Q\;\;$ if $\;\;0.97<Q\le 1.31$;} \\
    \hbox{$45.10+18.31Q\;\;$ if $\;\;1.31<Q<8.59$;} \\
    \hbox{$40.79+18.81Q\;\;$ if $\;\;Q\ge 8.59$.} \\
\end{array}%
\right.
$$

\section{Proof of the Theorems} \label{Proof}

We will consider the cases (i,j), where (i) is one of the cases of
Lemma \ref{L1}, and (j) is one of the cases of Lemma \ref{L2}.
First of all, suppose that at equilibrium the path $P_3$ has a
vanishing flow in $N^+$, i.e. the flow at equilibrium may only use
$P_1$ and $P_2$, or just one of those paths. It is easy to see
that the same flow is an equilibrium flow in $N$, and the travel
times at equilibria are equal. Thus, Braess' paradox cannot happen
in the following cases: (a,b), (a,c), (a,f), (b,b), (b,c), (b,f),
(c,b),(c,c),(c,f). Note that $T^+_{Eq} = T_{Eq}$ in the cases
(a,b), (b,c) and (c,f), while the conditions in other cases are
inconsistent. For example, in the case (a,c), Lemma \ref{L1} (a)
implies $\al<0$, while Lemma \ref{L2} (c) implies $\al>0$, so the
bounds are inconsistent.

Thus, there are 12 cases to consider. We will see that Braess'
paradox may occur in four cases only: (c,a), (c,d), (c,e) and
(c,g). This proves Mega-Theorem, while the cases themselves
provide proofs of Theorems \ref{T2}, \ref{T3}, \ref{T4} and
\ref{T1}, respectively.

\2 \noindent {\bf\em Case (a,a):} Assume that $T^+_{Eq} > T_{Eq}$,
which implies $Q(\be_{13}-\be_4)>\bar\al$. Since
$0<Q\le\bar\al/\be_{13}$, we obtain $\bar\al>0$. Therefore,
$\be_{13}-\be_4>0$ and
$$
Q > {\bar\al \over \be_{13}-\be_4} > {\bar\al \over \be_{13}},
$$
contrary to $Q\le\bar\al/\be_{13}$.



\2 \noindent {\bf\em Case (a,d):} $T^+_{Eq} > T_{Eq}$ implies
$\bar\al+Q\be_4<0$, i.e. $Q<-\bar\al/\be_4$, but
$Q>-\bar\al/\be_4$ in the case (d).

\2 \noindent {\bf\em Case (a,e):} We will show that the upper and
lower bounds imposed on $Q$ in this case are inconsistent.

Let us first assume that $\hat\al\ge0$. We have $Q\le
-\al/\be_{45}$ from (a) and $Q>\hat\al/\be_{35}$ from (e). Hence
$\hat\al/\be_{35} < -\al/\be_{45}$ or
 \beq \label{ae1}
\hat\al\be_{45} +\al\be_{35} <0.
 \eeq
From (e), we obtain
$$
{\hat\al \over \be_{35}} < {\bar{\al}\be_{25} + \al\be_{3} \over
\be_1\be_{25}+\be_3\be_{12}}.
$$
Using $\bar\al=\al+\hat\al$, the above inequality is equivalent to
$$
\hat\al(\be_1\be_{25}+\be_3\be_{12}-\be_{25}\be_{35})-\al\be_{35}\be_{235}<0.
$$
Since
$\be_1\be_{25}+\be_3\be_{12}-\be_{25}\be_{35}=(\be_1-\be_5)\be_{235}$,
we have
 \beq \label{ae2}
\hat\al(\be_{1}-\be_5) -\al\be_{35} <0.
 \eeq
Adding (\ref{ae1}) and (\ref{ae2}), we obtain $\hat\al\be_{14}<0$,
contrary to the assumption $\hat\al\ge0$.

Now suppose that $\hat\al<0$. We have $Q\le -\al/\be_{45}$ from
(a) and $Q>-\hat\al/\be_{2}$ from (e). Hence $-\hat\al/\be_{2} <
-\al/\be_{45}$ or
 \beq \label{ae3}
\hat\al\be_{45} -\al\be_{2} >0.
 \eeq
From (e), we obtain
$$
{-\hat\al \over \be_{2}} < {\bar{\al}\be_{25} + \al\be_{3} \over
\be_1\be_{25}+\be_3\be_{12}}.
$$
Hence, using $\bar\al=\al+\hat\al$,
$$
\hat\al(\be_1\be_{25}+\be_3\be_{12}+\be_{2}\be_{25})+\al\be_{2}\be_{235}>0.
$$
It is easy to see that
$\be_1\be_{25}+\be_3\be_{12}+\be_{2}\be_{25} = \be_{12}\be_{235}$,
so
 \beq \label{ae4}
\hat\al\be_{12} +\al\be_{2} >0.
 \eeq
Adding (\ref{ae3}) and (\ref{ae4}), we obtain $\hat\al\be>0$, a
contradiction.


\2 \noindent {\bf\em Case (a,g):} We will show that the upper and
lower bounds imposed on $Q$ in this case are inconsistent. We have
$\mu_1<Q\le-\al/\be_{45}$, i.e.
 \beq \label{ag1}
\hat\al\be_{45} +\al\be_{5} <0.
 \eeq
Also, from (g),
 \beq \label{QB}
Q{\cal B}_1 < \hat\al\be_{14} + \bar\al\be_{25}.
 \eeq

Suppose that ${\cal B}_1=0$, i.e.
 \beq \label{ag2}
\be_1\be_5=\be_2\be_4
 \eeq
and
$$
0<\hat\al\be_{14} + \bar\al\be_{25}.
$$
From (\ref{ag2}), we obtain $\be_1=\be_2\be_4/\be_5$,
$\be_2=\be_1\be_5/\be_4$ and $\be_1/\be_4=\be_2/\be_5$. Therefore,
$$
0<\hat\al\be_{14} + \bar\al\be_{25} =
\hat\al\be_{4}(1+{\be_2\over\be_5}) +
\bar\al\be_{5}(1+{\be_1\over\be_4}) = (\hat\al\be_{4}+
\bar\al\be_{5})(1+{\be_2\over\be_5}).
$$
Thus, $\hat\al\be_{4}+ \bar\al\be_{5}>0$. Now
$$0<\hat\al\be_{4}+
\bar\al\be_{5}= \hat\al\be_{4}+
(\al+\hat\al)\be_{5}=\hat\al\be_{45}+ \al\be_{5},
$$
 contrary to (\ref{ag1}).

Using $\bar\al=\al+\hat\al$, we can re-write (\ref{QB}) in the
following form:
$$
Q{\cal B}_1 < \al\be_{25} + \hat\al\be.
$$

Assume now that ${\cal B}_1<0$, i.e.
$$
Q > {\al\be_{25} + \hat\al\be\over {\cal B}_1}.
$$
Since $Q\le-\al/\be_{45}$, we obtain
$$
{\al\be_{25} + \hat\al\be\over {\cal B}_1} < {-\al\over \be_{45}}
$$
or
$$
\al({\cal B}_1+\be_{25}\be_{45}) + \hat\al\be_{45}\be >0.
$$
It is not difficult to see that the last inequality is equivalent
to
$$
\al\be_{5}\be + \hat\al\be_{45}\be >0
$$
or
$$
\al\be_{5} + \hat\al\be_{45} >0,
$$
contrary to (\ref{ag1}).

Finally, let us suppose that ${\cal B}_1>0$, i.e.
$$
Q < {\al\be_{25} + \hat\al\be\over{\cal B}_1}.
$$
Since $Q>\mu_1$, we obtain
$$
{\hat{\al}\be_{14} - \al\be_{3}\over L} < {\al\be_{25} +
\hat\al\be\over{\cal B}_1},
$$
where $L=\be_3\be_{45}+\be_5\be_{14}$. The last inequality is
equivalent to
$$
\hat\al(\be L-\be_{14}{\cal B}_1) + \al(\be_{25}L+\be_3{\cal B}_1)
>0.
$$
It is easy to see that
$$\be L-\be_{14}{\cal B}_1 =
\be_{45}(\be\be_3+\be_{14}\be_{25})$$
 and
$$
 \be_{25}L+\be_3{\cal
B}_1 = \be_{5}(\be\be_3+\be_{14}\be_{25}). $$
Thus,
$$
\hat\al\be_{45} +\al\be_{5} >0,
$$
contrary to (\ref{ag1}).


\2 \noindent {\bf\em Case (b,a):} Assume that $T^+_{Eq} > T_{Eq}$,
which implies $Q(\be_{35}-\be_2)>\hat\al$. Since
$0<Q\le\hat\al/\be_{35}$, we obtain $\hat\al>0$. Therefore,
$\be_{35}-\be_2>0$ and
$$
Q > {\hat\al \over \be_{35}-\be_2} > {\hat\al \over \be_{35}},
$$
contrary to $Q\le\hat\al/\be_{35}$.



\2 \noindent {\bf\em Case (b,d):} We will show that the upper and
lower bounds imposed on $Q$ in this case are inconsistent.

Let us first assume that $\bar\al\ge0$. We have $Q\le
\al/\be_{12}$ from (b) and $Q>\bar\al/\be_{13}$ from (d). Hence
$\bar\al/\be_{13} < \al/\be_{12}$ or
 \beq \label{bd1}
\bar\al\be_{12} - \al\be_{13} < 0.
 \eeq
From (d,) we obtain
$$
{\bar\al \over \be_{13}} < {\hat{\al}\be_{14} - \al\be_{3} \over
L},
$$
where $L=\be_3\be_{45}+\be_5\be_{14}$. Using
$\hat\al=\bar\al-\al$, the above inequality is equivalent to
$$
\bar\al(L-\be_{13}\be_{14})+\al\be_{13}\be_{134}<0
$$
or
 \beq \label{bd2}
\bar\al(\be_{5}-\be_1) + \al\be_{13} <0.
 \eeq
Adding (\ref{bd1}) and (\ref{bd2}), we obtain $\bar\al\be_{25}<0$,
contrary to the assumption $\bar\al\ge0$.

Now suppose that $\bar\al<0$. We have $-\bar\al/\be_{4} <
\al/\be_{12}$ or
 \beq \label{bd3}
\al\be_{4} + \bar\al\be_{12} >0.
 \eeq
From (d) we obtain
$$
{-\bar\al \over \be_{4}} < {\hat{\al}\be_{14} - \al\be_{3} \over
L},
$$
where $L=\be_3\be_{45}+\be_5\be_{14}$. Using
$\hat\al=\bar\al-\al$, the above inequality is equivalent to
$$
\bar\al(L+\be_{4}\be_{14})-\al\be_{4}\be_{134}>0
$$
or
 \beq \label{bd4}
\bar\al\be_{45} - \al\be_{4} > 0.
 \eeq
Adding (\ref{bd3}) and (\ref{bd4}), we obtain $\bar\al\be>0$,
contrary to the assumption $\bar\al<0$.

\2 \noindent {\bf\em Case (b,e):} $T^+_{Eq} > T_{Eq}$ implies
$\hat\al+Q\be_2<0$, i.e. $Q<-\hat\al/\be_2$, but
$Q>-\hat\al/\be_2$ in the case (e).


\2 \noindent {\bf\em Case (b,g):} We will show that the upper and
lower bounds imposed on $Q$ in this case are inconsistent. We have
$\mu_2<Q\le\al/\be_{12}$, i.e.
$$
\bar\al\be_{12} - \al\be_{1} <0,
$$
or, using $\al=\bar\al-\hat\al$,
 \beq \label{bg1}
\hat\al\be_{1} + \bar\al\be_{2} <0.
 \eeq
From (g), we have
 \beq \label{QB1}
Q{\cal B}_1 < \hat\al\be_{14} + \bar\al\be_{25}.
 \eeq

Suppose that ${\cal B}_1=0$, i.e.
 \beq \label{bg2}
\be_1\be_5=\be_2\be_4
 \eeq
and
$$
0<\hat\al\be_{14} + \bar\al\be_{25}.
$$
From (\ref{bg2}), we obtain $\be_4=\be_1\be_5/\be_2$,
$\be_5=\be_2\be_4/\be_1$ and $\be_4/\be_1=\be_5/\be_2$. Therefore,
$$
0<\hat\al\be_{14} + \bar\al\be_{25} =
\hat\al\be_{1}(1+{\be_5\over\be_2}) +
\bar\al\be_{2}(1+{\be_4\over\be_1}) = (\hat\al\be_{1}+
\bar\al\be_{2})(1+{\be_5\over\be_2}).
$$
Thus, $\hat\al\be_{1}+ \bar\al\be_{2}>0$,
 contrary to (\ref{bg1}).

Assume now that ${\cal B}_1<0$, i.e.
$$
Q > {\hat\al\be_{14} + \bar\al\be_{25}\over{\cal B}_1}.
$$
Since $Q\le\al/\be_{12}=(\bar\al-\hat\al)/\be_{12}$, we obtain
$$
(\bar\al-\hat\al) {\cal B}_1 <
\hat\al\be_{12}\be_{14}+\bar\al\be_{12}\be_{25}
$$
or
$$
\hat\al(\be_{12}\be_{14}+{\cal
B}_1)+\bar\al(\be_{12}\be_{25}-{\cal B}_1)
>0.
$$
The last inequality is equivalent to
$$
\hat\al\be_{1}\be+\bar\al\be_{2}\be >0,
$$
 contrary to (\ref{bg1}).

Finally, let us suppose that ${\cal B}_1>0$, i.e.
$$
Q < {\hat\al\be_{14} + \bar\al\be_{25}\over{\cal B}_1}.
$$
Since, $Q>\mu_2$, we obtain
$$
{\bar{\al}\be_{25} + \al\be_{3} \over \be_1\be_{25}+\be_3\be_{12}}
< {\hat\al\be_{14} + \bar\al\be_{25}\over{\cal B}_1}.
$$
Using $\al=\bar\al-\hat\al$, the last inequality can be re-written
as follows:
 \beq \label{bg3}
\hat\al X + \bar\al Y >0,
 \eeq
where
$$
X=\be_{14}(\be_1\be_{25}+\be_3\be_{12})+\be_3{\cal B}_1
$$
and
$$
Y=\be_{25}(\be_1\be_{25}+\be_3\be_{12}-{\cal B}_1)-\be_3{\cal
B}_1.
$$
It is not difficult to see that
$$
X= \be_1(\be_{14}\be_{25}+\be_3\be)
$$
and
$$
Y= \be_2(\be_{14}\be_{25}+\be_3\be).
$$
Thus, (\ref{bg3}) is equivalent to
$$
\hat\al \be_1 + \bar\al \be_2 >0,
$$
contrary to (\ref{bg1}).


\2 \noindent {\bf\em Case (c,a):}  This case corresponds to
Theorem \ref{T2}. It is easy to see that $T^+_{Eq} > T_{Eq}$ is
equivalent to
$$
Q(\be_{135}\be-\be_{12}\be_{45}) > \hat\al \be + \al \be_{12}
$$
or, using $\al=\bar\al-\hat\al$,
$$
Q{\cal B}_2 > \hat\al \be_{45} + \bar\al \be_{12}.
$$
Moreover,
$$
\max\left\{{\al\over\be_{12}}; {-\al\over\be_{45}}\right\} < Q \le
\min \left\{{\hat\al\over\be_{35}};
{\bar\al\over\be_{13}}\right\},
$$
which implies $\hat\al>0$ and $\bar\al>0$. Therefore, $\hat\al
\be_{45} + \bar\al \be_{12}>0$ and ${\cal B}_2 > 0$, i.e.
$$
Q> {\hat\al \be_{45} + \bar\al \be_{12}\over {\cal B}_2 }.
$$




\2 \noindent {\bf\em Case (c,d):}  This case corresponds to
Theorem \ref{T3}. It is easy to see that $T_{Eq} < T^+_{Eq}$ is
equivalent to
$$
\al_{12}+(\al+Q\be_{45})\be_{12}/\be <
\al_{45}+Q\be_{45}-(\bar\al+Q\be_4)\be_4/\be_{134}
$$
or
$$
Q{\cal B}_3 >  \bar\al\be_4 \be - \al\be_{134}\be_{45}.
$$
In addition, we have
$$
\max\left\{{\al\over\be_{12}}; {-\al\over\be_{45}}; {\bar\al\over
\be_{13}}; {-\bar\al\over\be_4} \right\} < Q \le \mu_1.
$$

Let us show that these inequalities are inconsistent if ${\cal
B}_3 \le 0$. Using (\ref{eq1}), the inequality
$-\al/\be_{45}<\mu_1$ can be written as follows:
$$
 {-\al\over\be_{45}} < {\bar\al \be_{14} -\al \be_{134}\over
 \be_3\be_{45}+\be_5\be_{14}},
$$
which is equivalent to
$$
\bar\al\be_{14}\be_{45} - \al (\be_{134}\be_{45}-
\be_3\be_{45}-\be_5\be_{14}) > 0
$$
or
$$
\bar\al > {\al \be_4 \over \be_{45}}.
$$
Therefore,
$$
Q{\cal B}_3 >  \bar\al\be_4 \be - \al\be_{134}\be_{45} > {\al
\be_4 \over \be_{45}}\be_4 \be - \al\be_{134}\be_{45} = {-\al
{\cal B}_3 \over \be_{45}}.
$$
Thus,
$$
Q{\cal B}_3 >   {-\al {\cal B}_3 \over \be_{45}},
$$
which is not satisfied if ${\cal B}_3 =0$. If ${\cal B}_3 <0$,
then we obtain
$$
Q <   {-\al \over \be_{45}},
$$
which is inconsistent with the inequality $Q > -\al /\be_{45}$. We
conclude that ${\cal B}_3 >0$ and hence
$$
Q > { \bar\al\be_4 \be - \al\be_{134}\be_{45}\over {\cal B}_3}.
$$

\2 \noindent {\bf\em Case (c,e):} This case corresponds to Theorem
\ref{T4}. It is not difficult to see that $T_{Eq} < T^+_{Eq}$ is
equivalent to
$$
\al_{12}+(\al+Q\be_{45})\be_{12}/\be <
\al_{12}+Q\be_{12}-(\hat\al+Q\be_2)\be_2/\be_{235}
$$
or
$$
Q{\cal B}_4 > \hat\al\be_2 \be + \al\be_{235}\be_{12}.
$$
Moreover, we obtain
$$
\max\left\{{\al\over\be_{12}}; {-\al\over\be_{45}}; {\hat\al\over
\be_{35}}; {-\hat\al\over\be_2} \right\} < Q \le \mu_2.
$$

Let us show that these inequalities are inconsistent if ${\cal
B}_4 \le 0$. Using $\bar\al=\al+\hat\al$, the inequality
$\al/\be_{12}<\mu_2$ can be written as follows:
$$
 {\al\over\be_{12}} < {\al \be_{235} +\hat\al \be_{25}\over
 \be_1\be_{25}+\be_3\be_{12}},
$$
which is equivalent to
$$
\al( \be_{12}\be_{235}-  \be_1\be_{25}-\be_3\be_{12}) + \hat\al
\be_{12}\be_{25} > 0
$$
or
$$
\hat\al > {-\al \be_2 \over \be_{12}}.
$$
Therefore,
$$
Q{\cal B}_4 >  \hat\al\be_2 \be + \al\be_{235}\be_{12} > {-\al
\be_2 \over \be_{12}}\be_2 \be + \al\be_{235}\be_{12} = {\al {\cal
B}_4 \over \be_{12}}.
$$
Thus,
$$
Q{\cal B}_4 >   {\al {\cal B}_4 \over \be_{12}},
$$
which is not satisfied if ${\cal B}_4 =0$. If ${\cal B}_4 <0$,
then
$$
Q <   {\al \over \be_{12}},
$$
which is inconsistent with the inequality $Q > \al /\be_{12}$.
Thus, ${\cal B}_4 >0$ and hence
$$
Q > {\hat\al\be_2 \be + \al\be_{235}\be_{12}\over {\cal B}_4}.
$$


\2 \noindent {\bf\em Case (c,g):} This case corresponds to Theorem
\ref{T1}. It is easy to see that $T^+_{Eq} > T_{Eq}$ is equivalent
to
$$
g{\cal B}_1/\be > 0.
$$
From the proof of Lemma \ref{L2} we know that $g>0$. Hence ${\cal
B}_1>0$. Therefore, the last inequality in the case (g) can be
written as
$$
 Q < {\hat\al\be_{14} + \bar\al\be_{25} \over {\cal
B}_1}.
$$
In addition, we have
$$
\max\left\{{\al\over\be_{12}}; {-\al\over\be_{45}}; \mu_1;
\mu_2\right\} < Q.
$$

\section{Conclusions and Further Research}

Braess' paradox has been investigated in symmetric four-link
networks by a number of authors. This paper provides a
generalisation of Pas and Principio's findings \cite{Pas} to
non-symmetric networks and shows, for a given congested network,
that Braess' paradox occurs if and only if the total demand of
travel ($Q$) lies between a certain range of values. The
motivation for such an extension was given by Pas and Principio in
their conclusions in \cite{Pas}. Also, in the context of
volume-delay functions and their parameters, it can be argued that
symmetry properties of networks  are not very common in real-life
situations.

A further research direction is to generalise results presented in
\cite{Pas} and this work, and find the necessary and sufficient
conditions for the existence of Braess' paradox with arbitrary
travel times of links in networks with complex topology. Further
extensions of this work might also include the investigation of
additional network features that were not illustrated in the
classical problem introduced by Braess and his colleagues
\cite{Bra}. Of specific relevance to transport researchers is the
study of networks in which volume-delay functions are not linear
in their parameters. Another future direction suggested by Pas and
Principio \cite{Pas} would be the investigation of Braess' paradox
in situations where the overall demand ($Q$) for travel is not
constant (i.e. inelastic).


\end{document}